\documentclass[11pt,twoside]{amsart}

\usepackage{amssymb}
\usepackage{amsmath}
\usepackage[latin1]{inputenc}

\theoremstyle{plain}

	\newtheorem{theorem}{Theorem}[section]
	\newtheorem{proposition}[theorem]{Proposition}
	{\unskip\nobreak\hskip 2em plus 1fil\nobreak%
  	\parfillskip=0pt \endtrivlist}

\newcommand{\CC}{\mathbb{C}}
\newcommand{\RR}{\mathbb{R}}
\newcommand{\Dirac}{/ \!\!\!\! D}
\newcommand{\Ah}{\operatorname{\widehat{A}}}
\newcommand{\clifford}{c}
\newcommand{\vol}{\mathrm{vol}}

\begin{document}

\title[Vafa--Witten bound on the complex projective space]{Vafa--Witten bound on the complex projective space}

\author{H\'el\`ene Davaux \& Maung Min-Oo}
\address{H\'el\`ene Davaux, Laboratoire UMPA, 46 all\'ee d'Italie, 69364 LYON, France}
\email{hdavaux@umpa.ens-lyon.fr}
\address{Maung Min-Oo, Department of Mathematics and Statistics, McMaster University, Hamilton, Ontario L8S 4K1, Canada}
\email{minoo@univmail.cis.mcmaster.ca}
\date{\today}

\begin{abstract}
We prove 
that the largest first eigenvalue
of the Dirac operator among all hermitian metrics 
on the complex projective space of odd dimension $m$, larger than 
the Fubini-Study metric is bounded by $\left(2m(m+1)\right)^{1/2}$.
\end{abstract}

\maketitle
%%%%%%%%%%%%%%%%%%%%%%%%%%%%%%%%%%%%%%%%%%%%%%%%%%%%%%%%%%%%%%%%%%%
\section{Introduction}
%%%%%%%%%%%%%%%%%%%%%%%%%%%%%%%%%%%%%%%%%%%%%%%%%%%%%%%%%%%%%%%%%%%

Lower and upper eigenvalue estimates for the Dirac operator on
a closed Riemannian spin manifold are derived by very different methods.
Lower bound estimate are usually based on a Bochner-Lichnerowicz-Weitzenb\"ock
formula and related to scalar curvature (\cite{Fri80, Kir90, Mor95, Sem99}).
In the other direction, C. Vafa and E. Witten \cite{VW84, Ati85} showed that there exists a common upper bound for the smallest eigenvalue of all {\sl twisted} Dirac 
operators on a given Riemannian manifold. The common upper bound they found 
depends on the metric of the base manifold but is independent of the bundle and the unitary connection that is used to twist the Dirac operator. 

The idea is as follows. Compare the Dirac operator $\mathcal D^0$ 
(or a multiple of it)
to a twisted Dirac operator $\mathcal D^1$ acting on sections of the same
vector bundles. By index theory make sure that $\mathcal D^1$ has a kernel.
Let $k$ the multiplicity of the eigenvalue $0$ of $\mathcal D^1$. Estimate
the norm of the difference (which is a zero order operator), $\|\mathcal D^0-
\mathcal D^1\| =\|L\|$, by geometric quantities. Then at least $k$ eigenvalues
of $\mathcal D^0$ are bounded by $\|L\|$.

Using this
method, H. Baum  \cite{Bau91} exhibited later an explicit upper bound for 
the first eigenvalue of the (untwisted) Dirac operator on an even-dimensional 
Riemannian manifold that can be sent on a sphere by a map of (high) non-zero 
degree. This bound depends on the Lipschitz norm of the map from the manifold 
to the round sphere. For the tori, there is also the paper of N. Anghel 
\cite{Ang00}.

In \cite{Her04}, Marc Herzlich gives an optimal result for the first
eigenvalue of the Dirac operator on the spheres.
The goal of this short note is to find an analogous result for the
complex projective spaces.

\begin{theorem}
\label{CPm}
Let $g$ be any hermitian metric and $f$ be the Fubini-Study metric 
of maximal sectional curvature $4$ on the complex projective space $\CC P^m$,
with $m$ odd.

If $g\geq f$, that is $g(v,v) \geq f(v,v)$ for all $v \in T \CC P^m$, then 
$$\lambda_1(g)^2 \leq 2 m (m+1)$$
where $\lambda_1(g)$ 
is the lower eigenvalue of the Dirac operator of $g$. 
\end{theorem}

There are also some techniques (the metric Lie derivative)
described by J-P. Bourguignon and P. Gauduchon \cite{BG92}
to compare spinors for different metrics. In particular, they prove
some results on eigenvalue of Dirac operator when the metric varies.
In the case of conformal change of the metric, we can give the
following corollary of
 O. Hijazi  work \cite{Hij86}: 

\begin{proposition} 
Let $(M,g_0)$ be a closed spin manifold. If $g = e^u g_0$ where
$u$ is a non negative function on $M$, then
$ \displaystyle |\lambda_1(g)| \leq |\lambda_1(g_0)|.$
\end{proposition}

\medskip
\noindent
We still do not know whether Theorem \ref{CPm} can be improved so that the result is optimal for the Fubini-Study metric.  That is, can we 
prove that, for $g \geq f$, we have
$$\lambda_1(g)^2 \leq (m+1)^2 = \lambda_1^2(f)$$
with equality if and only if $f=g$.\\

\noindent
{\small
\textit{Acknowledgement :}
The first author is grateful to Marc Herzlich for his advice and
encouragement}

%%%%%%%%%%%%%%%%%%%%%%%%%%%%%%%%%%%%%%%%%%%%%%%%%%%%%%%%%%%%%%%%%%%
\section{Background material}
%%%%%%%%%%%%%%%%%%%%%%%%%%%%%%%%%%%%%%%%%%%%%%%%%%%%%%%%%%%%%%%%%%%

%%%%%%%%%%%%%%%%%%%%%%%%%%%%%%%%%%%%%%%%%%%%%%%%%%%%%%%%%%%%%%%%%%%
\subsection{Tautological bundle}

Let $h$ be the standard hermitian product
on $\CC^{m+1}$. We denote 
$\langle \;,\; \rangle = \Re e h.$ It is the standard riemannian product
on $\RR^{2m+2} = \CC^{m+1}$.
 
If $\CC P^m = \{ [z], z \in S^{2m+1} \subset \CC^{m+1}\}$, then
the tautological bundle is given by
$$E = \{([z],w) , \;\; w \in [z]\} \subset \CC P^{m} \times \CC^{m+1}.$$
We consider the complement $F$ of $E$:
$$F = \{([z],w) , \;\; w \perp_h z\} \subset \CC P^{m} \times \CC^{m+1}.$$
We have $E \oplus F = \CC P^{m} \times \CC^{m+1}.$
Coming from $\CC^{m+1}$, there are an hermitian structure
on the bundles $E$ and $F$.

We denote by $\alpha$ the standard identification between
$T \CC P^m$ and $\mathrm{Hom}_\CC(E,F)$ defined as follow:
%$$\alpha : T \CC P^m \cong \mathrm{Hom}_\CC(E,F)$$
Let $X$ be a vector of $T_x\CC P^m$,  
$X$ is generated by the curve $x(t) \in \CC P^m$ such 
that $x(0)=x$ and $\dot x (0)=X$. We associate to $X \in T_x \CC P^m$ the 
homomorphism 
$$\begin{array}{rrcl}
\alpha(X) : & 
E_x \cong x &\rightarrow & x^\perp = F_x\\
  & u &\mapsto &  \Pi_F( \nabla_X^0(u)) \quad \textrm{ where } \nabla_X^0(u) =\dot u(0)
\end{array}
$$ 
where
$u(t)$ is a curve of $\CC^{m+1}-\{0\}$, $u(0) =u$ and 
$\forall t, u(t) \in x(t)$. For another curve $v$, we check that
$\dot u(0)=\dot v(0)$ in $\CC^{m+1}/x$. We identify $\CC^{m+1}/x$
with $x^\perp$ where we use the standard hermitian product $h$
on $\CC^{m+1}$ to define the orthogonal.

Using this identification $\alpha$, 

$\bullet$ we can give a definition for the Fubiny-Study
metric. Let $X,Y \in T_x\CC P^m$,
$$f_x(X,Y)= \Re e \left[h(\alpha(X)(u), \alpha(Y)(u))\right] = 
\langle \alpha(X)(u), \alpha(Y)(u) \rangle$$
where $u$ is a unitary generator of $x$
(the expression of this $f$ in local coordinate coincides with the
usual one described in \cite[T. 2, chap. IX, \S 6,7]{KN69}).

$\bullet$ let $\epsilon$ be a section of $E$, 
the connection for the bundle $E$ can be written as:
$$\nabla^E_V \epsilon = \Pi_E(\nabla^0_V\epsilon) =
\nabla^0_V\epsilon - \alpha(V)(\epsilon)$$
where $\nabla^0$ is the flat connection in the
trivial bundle $\CC P^m \times \CC^{m+1}$
(so $\nabla^0_V \epsilon$ is the usual derivative $V(\epsilon)$ of 
the function
$\epsilon : \CC P^m \rightarrow \CC^{m+1}$).

%%%%%%%%%%%%%%%%%%%%%%%%%%%%%%%%%%%%%%%%%%%%%%%%%%%%%%%%%%%%%%%%%%%
\subsection{Hilbert polynomial}

Let consider $H$ the dual of $E$. Let $\gamma$ be a $2$-form on $\CC P^m$ 
such that $[\gamma] = c_1(H)$.  
It is well known that  $\int_{\CC P^1} \gamma = -1$.
We can compute $c_1(M) = -(m+1) [\gamma]$ so $\CC P^m$ is spin iff $m$ is odd.
We have $c_1(\kappa) = -c_1(M)$ where 
$\kappa$ is the canonical bundle (i.e. $\kappa = \wedge^n_\CC T^*M$).
 If $M$ is spin, $c_1(\kappa)\equiv 
0 [2]$ thus there exists a square root for the line bundle $\kappa$, which
we denote $\kappa^{1/2}$. As $c_1(H)=[\gamma]$, 
$c_1(\kappa)= [(m+1)\gamma]$ and 
 the line bundle are classified
by there first Chern class, we have $H^{\otimes (m+1)} \cong \kappa$ and 
$H^{\otimes (m+1)/2} \cong \kappa^{1/2}$.

We define and compute  the Hilbert polynomial (\cite[IV \S 11, p. 365]{LM89})
$$P(t) := \{e^{\frac 1 2 t \gamma} \Ah(M)\}[M] = 
\frac{1}{2^m m!} \prod_{j=1}^{m} (t -m+2j-1)$$
where $\Ah(M)$ is the $\Ah$-class of $M$ belonging to the cohomology
$H^{4*}(M)$.

For any metric $g$ on $\CC P^m$, we let $\Sigma_g$ be its spin bundle
and $\Dirac^g$ is the Dirac operator relative to the metric $g$. 
We apply Atitah-Singer index theorem to the twisted Dirac operator 
%\begin{equation*}
$\Dirac^{g+}_{\kappa^{1/2}}$. 
%: \Sigma^+_g \otimes\kappa^{1/2}
% \longrightarrow \Sigma^-_g \otimes \kappa^{1/2}.
%\end{equation*}
One of the most important consequences of index theory is 
the topological invariance of the index which is given by
$$\textrm{ind}\left(  \Dirac^{g+}_{\kappa^{1/2}} : 
\Sigma_g^+\otimes  \kappa^{1/2}\longmapsto \Sigma_g^-\otimes  \kappa^{1/2}
\right) = \{\textrm{ch}(\kappa^{1/2}) \Ah(M)\}[M]= P(m+1)=1$$
where  $\textrm{ch}(\kappa^{1/2})= e^{c_1(\kappa^{1/2})}$ is the Chern 
character of the line bundle $\kappa^{1/2}$.

We have also
$P(-m-1) = -1$. Therefore the twisted Dirac operator on 
$\Sigma_g \otimes E^{\otimes (m+1)/2}$ has non-zero harmonic spinors.

%%%%%%%%%%%%%%%%%%%%%%%%%%%%%%%%%%%%%%%%%%%%%%%%%%%%%%%%%%%%%%%%%%%
\subsection{Complementary bundle for $E^{\otimes (m+1)/2}$}
 
We remark that, 
given a classifying map $\phi : X \rightarrow \CC P^m$ for a line 
$L = \phi^* E$, we can defined
 $L^\perp:= \phi^* F$, then  the bundle $L \oplus L^\perp$ is trivial,
that is $L^\perp$ is a complementary bundle for $L$.

Let $k = \frac{m+1}{2}$. 
Consider all the homogeneous monomials of degree $k$ 
with $(m+1)$ variables $X_0,X_1, \cdots X_m$ and more precisely the 
folowing monomials:
$$M_\alpha(X_0,X_1, \cdots, X_m) = 
\sqrt{\frac{k!}{\alpha_0! \cdots \alpha_m!}} 
X_0^{\alpha_0}X_1^{\alpha_1} \cdots X_{m}^{\alpha_m}$$ 
where $\sum \alpha_i =k$. There are 
$N+1$ such monomials with
$ \displaystyle N+1 =\left( \begin{array}{c} m+k \\ m \end{array}\right).$

Now we define the Veronese map of degree $k$ given in
homogeneous coordinates by
$$
\begin{array}{rrcl}
\psi_k  : & \CC P^m &\longrightarrow & \CC P^N\\
&[z_0 : \cdots : z_m] &\longmapsto & [M_0(z_0, \cdots , z_m) :
 \cdots : M_N(z_0, \cdots , z_m)]
\end{array} 
$$
The map $\psi_k$  
is an embedding of $\CC P^m$ in $\CC P^N$. 
Thanks to the coefficient $\sqrt{\frac{k!}{\alpha_0! \cdots \alpha_m!}}$
(see  \cite{MS82,NO76}), 
we can lift this map to the spheres and, by a computation in local coordinate,
we can prove that
$$\psi_k^*f_{\CC P^N} = k f_{\CC P^m}.$$

Moreover, the map $\psi_k$  classifies 
$H^{\otimes k} = \psi_k^*(H)$ on $\CC P^m$
(see \cite[chap. 1, \S4]{GH78}).
 Therefore
$$E_k:= \psi_k^*(E) \cong E^{\otimes k}, \;\;
F_k := \psi_k^*(F)= E^{\otimes k \perp} \;\; \textrm{ and } \;\; 
E_k \oplus F_k = \CC P^m \times \CC^{N+1} =:\mathbb T$$
We have also the pull-back 
connections 
$\nabla^{E_k}:=\psi_k^*(\nabla^E)$ on $E_k$ and 
$\nabla^{F_k}:=\psi_k^*(\nabla^F)$ on $F_k$.
Thus, on $\mathbb T = E_k \oplus F_k$, we have two connections 
$\psi_k^* \nabla^0$ (it is the trivial connection)
and $\psi_k^* (\nabla^E \oplus \nabla^F)$.

%%%%%%%%%%%%%%%%%%%%%%%%%%%%%%%%%%%%%%%%%%%%%%%%%%%%%%%%%%%%%%%%%%%
\section{The proof}
%%%%%%%%%%%%%%%%%%%%%%%%%%%%%%%%%%%%%%%%%%%%%%%%%%%%%%%%%%%%%%%%%%%

%%%%%%%%%%%%%%%%%%%%%%%%%%%%%%%%%%%%%%%%%%%%%%%%%%%%%%%%%%%%%%%%%%%
\subsection{Review on the approach of Vafa-Witten}

As $E_k \cong E^{\otimes (m+1)/2}$, 
we have already  proved, 
using Hilbert polynomial, that the index 
of 
$\Dirac^{g+}_{E_k}$
is non zero and
thus the  kernel of $\Dirac^{g}_{E_k}$ is non trivial.
Let us now consider the tensor product bundle 
$\Sigma_g\otimes (E_k \oplus F_k)$
endowed
with the connection $\nabla^g\otimes 1 + 1\otimes (\nabla^{E_k} \oplus 
\nabla^{F_k})$ and denote by $\mathcal D^1$ 
the twisted Dirac operator  
attached to this connection. Its kernel is also non trivial
and there exists a non zero harmonic spinor $\phi_0$ of $\mathcal D^1$
which lives in the $L^2$-sections of $\Sigma_g\otimes E_k \subset 
\Sigma_g\otimes (E_k \oplus F_k) = \Sigma_g\otimes \mathbb T$.

Let us now consider the tensor product bundle 
$\Sigma_g\otimes \mathbb T$, 
endowed
this time with the connection $\nabla^g\otimes 1 + 1\otimes 
\psi_k^*\nabla^0$. 
As the pair 
$(\mathbb T,\psi_k^*\nabla^0)$ is a \emph{trivial} flat bundle on 
$\CC P^m$, the spectrum of 
the twisted Dirac operator $\mathcal D^0$ 
attached to this connection is the same as 
the spectrum 
of the Dirac operator $\Dirac^g$ on $\Sigma_g$, 
but with each eigenvalue repeated $N+1$ times 
its multiplicity. 

Let consider the zero order operator 
$ \displaystyle  L := \mathcal{D}^{0} - 
\mathcal{D}^{1}.$
For a non-zero harmonic spinor $\phi_0$ of $\mathcal D^1$, which
belongs to the $L^2$-sections of
 $\Sigma_g \otimes \psi_k^*E$, we have
%$$
\begin{eqnarray*}
|\lambda_1(g)| 
\|\phi_0\|^2 \leq |\langle \mathcal D^0 \phi_0 , \phi_0 \rangle|
= |\langle \mathcal{D}^{0}\phi_0 - \mathcal{D}^{1} \phi_0, \phi_0 \rangle|
\leq \|L\| \|\phi_0\|^2
\end{eqnarray*}
%$$
where $L$ is considered as a linear operator from the 
$L^2$-sections of $\Sigma\otimes \psi_k^*E$ to the 
$L^2$-sections of $\Sigma \otimes \mathbb T$.
Thus
$$\lambda_1(g)^2\leq \|L\|^2.$$

%%%%%%%%%%%%%%%%%%%%%%%%%%%%%%%%%%%%%%%%%%%%%%%%%%%%%%%%%%%%%%%%%%%
\subsection{Computation}
We have
$$\|L\|^2 = \sup_{\|\phi\|=1} \|L(\phi)\|^2
=\sup_{\|\phi\|=1} \int_{\CC P^m} \langle L_x (\phi_x), L_x (\phi_x) \rangle
\mathrm d \vol_g \leq \max_{x\in \CC P^m} \|L_x\|^2.$$
Thus is suffices to give an upper bound for $\|L_x\|$: the computation
is pointwise.

It  is easy to  compute $L_x$ on a decomposed section $\sigma\otimes 
\psi_k^*\tau \in
\Sigma_g \otimes \psi_k^*E \subset \Sigma_g \otimes \mathbb T$. Moreover,
it is sufficient because  $E$ is of complex rank $1$ : in fact, 
we can write every complex section of $\Sigma_g \otimes E_k$ as
 $\sigma \otimes \psi_k^*\tau$ where $\sigma$ is a section of $\Sigma_g$
and $\tau$ is a fixed everywhere non-zero section of $E$.

We denote the Clifford action relative to $g$ on the bundle
$\Sigma_g$ by $\clifford_g(\cdot)$.
If $\{e_a\}_{a= 1, \cdots, 2m}$ is any $g$-orthonormal basis on $\CC P^m$,
\begin{eqnarray*}
L_x (\sigma\otimes\psi_k^*\tau) & 
= &\left[ \sum_{a=1}^{2m} \, \clifford_g(e_a) \big(
\nabla^g_{e_a}\sigma\otimes\psi_k^*\tau + 
\sigma\otimes(\psi_k^*\nabla^0)_{e_a}\psi_k^*\tau \right.\\ 
&& \left. \hspace{2cm} - 
\nabla^g_{e_a}\sigma\otimes\psi_k^*\tau
 - \sigma\otimes (\psi_k^*\nabla^E)_{e_a}\psi_k^*\tau \big)
 \right]_x\\
& = & \sum_{a=1}^{2m} \, \left[\clifford_g(e_a)  
\sigma\right]_x \otimes\left[\nabla^0_{\psi_{k*}e_a}\tau 
 - \nabla^E_{\psi_{k*}e_a}\tau \right]_{\psi_k (x)}\\
& = & \, \sum_{a=1}^{2m} \, \left[\clifford_g(e_a)
\sigma\right]_x \otimes 
\left[\alpha(\psi_{k*}e_a) (\tau)\right]_{\psi_k(x)}
\end{eqnarray*}

As $g \geq f$,
we can choose a $f$-orthonormal basis $\{\varepsilon_a\}$ 
and a $g$-orthonormal basis $\{e_a\}$ such that $e_a = \mu_a \varepsilon_a$
with $\mu_a \in ]0,1]$. Moreover, as the metric $g$ is hermitian,
we can suppose that $Je_a = e_{m+a}$ for $a \in \{1,\cdots ,m\}$. Thus
$\mu_{m+a}=\mu_a$.

We use the hermitian metric induced by
$g$ on $\Sigma_g$ and  the hermitian metric $h$ coming from
$\CC^{N+1}$ on $\mathbb T$.  We have:
\begin{eqnarray*}
\|L_x(\sigma \otimes \psi^*_k\tau)\|^2 & =
\sum_{a,b =1}^{2m} 
\langle \clifford_g(e_a)\sigma, \clifford_g(e_b) \sigma \rangle_{g,x} \cdot
h( \alpha({\psi_k}_*(e_a))(\tau),\alpha({\psi_k}_*(e_b))(\tau))_{\psi_k(x)}
\end{eqnarray*}

A short computation shows that
$\langle \clifford_g(e_a)\sigma, \clifford_g(e_a) 
\sigma \rangle_{g,x}= |\sigma|^2_g$
and $\langle \clifford_g(e_a)\sigma, \clifford_g(e_b) 
\sigma \rangle_{g,x}$ is a pure imaginary complex if  $a \neq b$,
thus 
\begin{eqnarray*}
\|L_x(\sigma \otimes \psi^*_k\tau)\|^2& =&
\sum_{a =1}^{2m} |\sigma|^2_{g,x}
 \cdot
h( \alpha({\psi_k}_*(e_a))(\tau),\alpha({\psi_k}_*(e_a))(\tau))_{\psi_k(x)}\\
&&+\sum_{a\neq b} 
\langle \clifford_g(e_a)\sigma, \clifford_g(e_b) 
\sigma \rangle_{g,x} \cdot i \Im m 
\left[h( \alpha({\psi_k}_*(e_a))(\tau),\alpha({\psi_k}_*(e_b))(\tau))_{\psi_k(x)}\right]
\end{eqnarray*} 

Using the definition of the Fubini-Study metric in terms of $\alpha$,
the equality $\psi_k^*f^{\CC P^N}=k f^{\CC P^m}$ and  the fact that
$\{\frac{e_a}{\mu_a}\}=
\{\epsilon_a\}$ is an orthonormal basis for $f^{\CC P^m}$, we prove that:
\begin{eqnarray*}
h( \alpha({\psi_k}_*(e_a))(\tau),\alpha({\psi_k}_*(e_a))(\tau))_{\psi_k(x)}
&=& \mu_a^2 k |\tau|^2 
\end{eqnarray*}
and
\begin{eqnarray*}
\Im m h
( \alpha({\psi_k}_*(e_a))(\tau),\alpha({\psi_k}_*(e_b))(\tau))_{\psi_k(x)}
&=&\left\{ \begin{array}{l}
-\mu_a^2 k |\tau|^2  \textrm{ if } a \in \{1, \cdots,m\} \textrm{ and } b=m+a\\
+\mu_b^2 k |\tau|^2  \textrm{ if } b \in \{1, \cdots,m\} \textrm{ and } a=m+b \\
0 \textrm{ in the other cases}
\end{array}
\right.
\end{eqnarray*}

Thus, using $\mu_a = \mu_{m+a}$,
\begin{eqnarray*}
\|L_x(\sigma \otimes \psi^*_k\tau)\|^2 & =&
2k |\sigma|^2 |\tau|^2 \sum_{a =1}^{m} \mu_a^2  \\
&+&k |\tau|^2 \left[-\sum_{a=1}^m i \mu_a^2 
\langle \clifford_g(e_a)\sigma, \clifford_g(Je_a) 
\sigma \rangle_{g,x}
+ \sum_{b=1}^m i \mu_b^2 
\langle \clifford_g(Je_b)\sigma, \clifford_g(e_b) 
\sigma \rangle_{g,x} \right]\\
&=& k |\tau|^2 \left[ 
2 \sum_{a =1}^{m} \mu_a^2 |\sigma|^2 
+ 2i\langle \clifford_g(\Omega) \sigma,\sigma \rangle
\right]
\end{eqnarray*}
where $\Omega$ is the K\"ahler form relative to $f$ and
$\clifford_g(\Omega) := -\sum_{a=1}^m \mu_a^2  \clifford_g(e_a)
\clifford_g(Je_a)$.

We check that
$$\|\clifford_g (\Omega)^2 \| \leq (\sum_{a=1}^m \mu_a^2)^2,
\quad \textrm{ thus } \quad
\|\clifford_g (\Omega) \| \leq \sum_{a=1}^m \mu_a^2.$$
As $\mu_a\leq 1$ and $k = (m+1)/2$, we conclude that
$$ \|L\|^2 \leq 4k\sum_{a=1}^m \mu_a^2 \leq  2m(m+1).$$

\end{document}